\documentclass[12pt]{article} 
\usepackage{amsfonts}

\usepackage[utf8]{inputenc} 

\usepackage{geometry} 
\usepackage{graphicx} 
\usepackage{xcolor}

\usepackage{booktabs} 
\usepackage{array} 
\usepackage{paralist} 
\usepackage{verbatim} 
\usepackage{subfig} 
\usepackage{biblatex}
\usepackage{hyperref}

\usepackage[english]{babel}
\usepackage[title]{appendix}%

\usepackage{amsmath}
\usepackage[utf8]{inputenc}
\usepackage{authblk}
\usepackage{fancyhdr} 
\usepackage{sectsty}

\usepackage[nottoc,notlof,notlot]{tocbibind} 
\usepackage[titles,subfigure]{tocloft} 
\usepackage{indentfirst}

\usepackage{amsthm}

\theoremstyle{plain}
\newtheorem{theorem}{Theorem}
\newtheorem{lemma}{Lemma}[section]

\newtheorem{proposition}{Proposition}[section]
\newtheorem{corollary}{Corollary}[section]
\theoremstyle{definition}

\theoremstyle{remark}

\title{Cover-time Gumbel Fluctuations in Finite-Range, Symmetric, Irreducible Random Walks on Torus}
\author{Hao Ge$^{1}$, Xiao Han$^2$, Yuan Zhang$^3$}
\date{ $^1$ Beijing International Center for Mathematical Research(BICMR) and Biomedical Pioneering Innovation Center(BIOPIC), Peking University, Beijing 100871, China\\
    $^2$Section de mathématiques, École Polytechnique Fédérale de Lausanne, CH-1015 Lausanne, Switzerland\\
    $^3$ Center for Applied Statistics and School of Statistics, Renmin University of China, Beijing 100872, China} 

\begin{document}
\maketitle
\abstract{In this paper, we rigorously establish the Gumbel-distributed fluctuations of the cover time, normalized by the mean first passage time, for finite-range, symmetric, irreducible random walks on a torus of dimension three or higher. This has been numerically demonstrated in (Chupeau et al. Nature Physics, 2015), supporting the broader applicability of the Gumbel approximation across a wide range of stochastic processes. Expanding upon the pioneering work of Belius (Probability Theory and Related Fields, 2013) on the cover time for simple random walks, we extend the proof strategy to encompass more general scenarios. Our approach relies on a strong coupling between the random walk and the corresponding random interlacements. The presented results contribute to a better understanding of the cover-time behavior in random search processes.
}

\section{Introduction}
\subsection{Cover Time and Gumbel Fluctuations}

Random search is a simple and straightforward optimization technique used to find optimal or approximate solutions, and is often employed in machine learning, optimization problems and parameter tuning tasks. The idea behind random search is to explore the search space by sampling points randomly, rather than following a specific systematic or deterministic pattern \cite{15,16,17}.

 Cover time, on the other hand, is a concept used in the study of random walks \cite{18,19,20}. The cover time represents the number of steps required for a random walk to visit all the vertices in a given graph or network. When using random search to explore a search space, we can view the process as a random walk through the space, where each point sampled represents a step in the walk. The cover time then corresponds to the number of iterations or samples required for the random search process to cover or explore the entire search space. The study of cover time helps in understanding the efficiency and convergence properties of random search algorithms and provides insights into their performance characteristics.

The cover time is defined as the maximal hitting time among all sites in the underlying set. Namely, for a random walk $X_t$, $t\geq 0$ on a finite set $S$, the cover time of it is defined by:
$$
T_C:=\max_{x\in S} H_x,
$$
where $H_x$ is the hitting time of $x$. Furthermore, for any $F \subset S$,
$$
T_C^F:=\max_{x\in F} H_x.
$$
Extensive research has been conducted on the cover time of various types of random walks, revealing the presence of Gumbel fluctuations. We say a sequence of random walks on (possibly different) finite sets exhibits Gumbel fluctuations if there exist two sequences of constants $\{a\}$ and $\{b\}$ such that the normalized cover time $\frac{T_C-a}{b}$ converges in distribution to the standard Gumbel distribution or extreme value distribution $G$, with the cumulative distribution function $F(z)=e^{-e^{-z}}$, as the size of state space goes to infinity. This phenomenon, akin to a central limit result, has been observed in various scenarios, including the well-known ``coupon collector's problem" which can be equivalently seen as the cover time problem for a random walk with laziness on a complete graph and more general cases, see e.g. Example 3.6.11 in \cite{13} and Theorem 7.9 in \cite{11}.

The study of cover time for simple random walks on discrete torus of dimension three or higher has long been a subject of interest. David Belius proved this Gumbel fluctuation conjecture \cite{1}, using the method of random interlacements introduced by Sznitman in 2011 to analyze the percolative properties of simple random walks \cite{5}. Furthermore, numerical experiments conducted by Chupeau et al. in 2015 provided additional evidence supporting the prevalence of Gumbel fluctuations in various types of random walks on tori \cite{2}. It suggests that such Gumbel fluctuation seems to widely exist in different types of random walk on tori.

\subsection{Our main results}

In this paper, we establish the Gumbel fluctuations of the cover time of finite range, symmetric, irreducible random walks in a discrete torus $\mathbb{T}_N$ with dimension at least three and length $N>0$. This finding provides a comprehensive explanation for the examples presented in \cite{2}, supporting the broader applicability of the Gumbel approximation across a wide range of stochastic processes. Our result also serves as a generalization of \cite{1}, and we draw inspiration from their proof technique to establish our own.

We now introduce the notations for a finite range, symmetric, irreducible random walk. One may see \cite{4} for more details. We call $V = \{x_1, . . . , x_l\} \subset \mathbb{Z}^d \setminus \{0\}$ a finite generating set if each $y\in \mathbb{Z}^d$ can be written as $k_1x_1 + ... + k_lx_l$ for some $k_1, k_2, ..., k_l \in \mathbb{Z}$. We let $\mathcal{G}$ denote the collection of all finite generating sets $V$ with the property that for all $x = (x^{(1)}, . . . , x^{(d)}) \in V$, the first nonzero component of x is positive. An example of such a set is the canonical orthogonal basis $\{e_1, . . . , e_d\}$ where for any $i,j \in \{1, 2, ..., d\}$,
\begin{equation}
\label{ei}
e_i^{(j)}=1_{i=j}
\end{equation}

Note that a finite range, symmetric, irreducible random walk in $\mathbb{Z}^d$ is given by specifying a 
$V =\{x_1, . . . , x_l\} \in \mathcal{G}$ and a function $\kappa : V \to (0, +\infty)$. To show this we further define the symmetric transition rate matrix, or the moving mode matrix $Q=(q_{ij})_{i,j \in \mathbb{Z}^d}$, to be
\begin{equation}
\label{1.1}
q_{ij}=\left\{
\begin{aligned}
\frac{1}{2} \kappa(i-j) & ,  &\indent  i-j \in V , \\
\frac{1}{2} \kappa(j-i) & ,  &\indent  j-i \in V ,\\
-\sum_{x\in V}\kappa(x) & , &\indent i=j ,\\
0 & , &\indent \text{otherwise}.\\
\end{aligned}
\right.
\end{equation}
The continuous-time Markov jumping process $Z_t$ induced by $Q$ is actually the finite range, symmetric, irreducible random walk in $\mathbb{Z}^d$ induced by $Q$. If we define $L:\mathbb{Z}^d \to \mathbb{T}_N$ by $L(x^{(1)},x^{(2)}, ... ,x^{(d)}) = (x^{(1)}+N\mathbb{Z},x^{(2)}+N\mathbb{Z},..,x^{(d)}+N\mathbb{Z})$, we further get all the finite range, symmetric, irreducible random walk in $\mathbb{T}_N$ by writing $Y_t :=L(Z_t)$. We also call $Q$ a moving mode matrix of $Y_t$. Let $P_x$ be the law of $Y_t$ starting from $x \in \mathbb{T}_N$, $P$ be the law of $Y_t$ starting from uniform distribution in $\mathbb{T}_N$ and $P_x^{\mathbb{Z}^d}$ be the law of $Z_t$ starting from $x \in \mathbb{Z}^d$. To simplify the problem we always assume $\sum_{x \in V}{\kappa(x)}=1$. 

Our main result is that

\begin{theorem}
\label{T1}
For any $d \geq 3$, given the moving mode matrix $Q$, there exist constants $c, c' >0$ such that for any $N>0$, $F \subset \mathbb{T}_N$ we have
\begin{equation}
\label{2.22}
\sup_{z \in \mathbb{R}}|P[T_C^F \leq N^d u_F(z)] -e^{-e^{-z}}| \leq c' |F|^{-c}.
\end{equation}
\end{theorem}
\noindent See (\ref{1.21}) for the definition of $u_F(z)$. Note that if we take $F=\mathbb{T}_N$ and $N\to \infty$, (\ref{2.22}) straightforwardly implies that
\begin{corollary}
For any $d \geq 3$, given the moving mode matrix $Q$, when $N\to \infty$ we have
\begin{equation}
\label{0.2}
\frac{T_C}{g(0)N^d}-\log N^d \stackrel{d}{\rightarrow} G.
\end{equation}  
\end{corollary}

Furthermore, we note that Proposition 3.7 in \cite{14} implies $g(0)N^d$ and the mean first passage time have the same order for simple random walk in $\mathbb{T}_N$ when $N\to \infty$. This result as well as its proof can be quite straightforwardly generalized into finite-range symmetric irreducible random walks. Namely, taking $V$ there to be a one-point set $\{x\}\subset \mathbb{T}_N$, one can get
\begin{proposition}
For any $d \geq 3$, given the moving mode matrix $Q$, when $N\to \infty$ we have
 \begin{equation}
\label{0.3}
\frac{\mathbb{E}H_x}{g(0)N^d}\to 1,
\end{equation}
\end{proposition}

 Thus combining (\ref{0.2}) and (\ref{0.3}), we further arrive at the following corollary predicted by the numerical simulations in \cite{2}.
 \begin{corollary}
 For any $d \geq 3$, given the moving mode matrix $Q$, when $N\to \infty$ we have
 \begin{equation}
 \frac{T_C}{\mathbb{E}H_x}-\log N^d \stackrel{d}{\rightarrow} G.
 \end{equation}
 \end{corollary}

Following a similar strategy as in \cite{1}, we establish a crucial coupling between random walks on the torus with the same moving mode matrix and the corresponding random interlacements for general weighted graphs defined in \cite{7}. This coupling has independent significance and is of particular interest. 

Recently, Berestycki et al. in \cite{21} proved under a broader setting the Gumbel fluctuations of random walks on general finite vertex-transitive graphs using a different approach relying on a finitary version of Gromov's theorem \cite{23}. Our work is independent of theirs and is more of a natural generalization of \cite{1} alongside the random interlacements approach. The coupling method in our approach gives better estimate in the setting of finite range, symmetric, irreducible random walks on torus.

The remainder of the paper is organized as follows. In Section 2, we introduce the necessary notations and preliminary lemmas. Section 3 presents the key coupling result together with its proof. In Section 4, we utilize the coupling to conclude the cover time estimate.

Throughout this paper, constants denoted by $c$, $c'$, etc., may vary in different contexts while constants represented by $c_1, c_2, \ldots$ have fixed values. All constants are positive and, unless otherwise specified, depend only on the dimension $d$ and the moving mode matrix $Q$ (see (\ref{1.1}) for its definition). Any constant that depends on a parameter, such as $\alpha$, is denoted as $c(\alpha)$.

\section{Preliminaries}

We introduce some basic notations and lemmas in this section.

We denote the $d$-dimensional discrete torus of side length $N \geq 3$ by $\mathbb{T}_N=(\mathbb{Z}/N\mathbb{Z})^d$. For two sets $W_1,W_2 \subset \mathbb{T}_N$ , we denote their difference as $W_1 \setminus W_2:=\{y: y \in W_1, y \notin W_2\}$. For each point $x \in \mathbb{T}_N$ and any set $W\subset \mathbb{T}_N$, we denote spatial shifts as follows: $W-x:=\{y-x:y\in W\}$, $W+x:=\{y+x:y\in W\}$.
For $x \in \mathbb{Z}^d$, we write $|x|_{\infty}$ for the $l_{\infty}$ norm of x. 
We use $d_{\infty}(\cdot,\cdot)$ to denote the distance 
in $\mathbb{T}_N$ induced by $|\cdot |_{\infty}$ and the periodic boundary condition. The closed $l_\infty$-ball of radius $r \geq 0$ with center $x$ in $\mathbb{Z}^d$
 or  $\mathbb{T}_N$ is denoted by $B(x,r)$. 

For a set $U$ ($U=\mathbb{Z}^d$ or $\mathbb{T}_N$) we write $\Gamma(U)$ for the space of all cadlag piecewise constant functions from $[0,\infty)$
to $U$, with at most a finite number of jumps in any finite time interval. For convenience we sometimes write $Y(a, b)$ to denote the set $\{Y(t) : t \in [a, b] \cap [0,\infty)\}$ for a trajectory $Y  \in \Gamma(U)$.

For a finite measure $\nu$ on $\mathbb{T}_N$ or $\mathbb{Z}^d$, we sometimes write $P_\nu:=\sum_{x\in\mathbb{T}_N} \nu(x) P_x$,  $P_\nu^{\mathbb{Z}^d}:=\sum_{x\in\mathbb{Z}^d} \nu(x) P_x^{\mathbb{Z}^d}$. We also define the entrance time, return time and exit time of $Y_t$ by
\begin{equation}
\label{1.2}
H_A=\inf \{ t\geq 0:Y_t \in A\}, \tilde{H}_A=\inf \{ t\geq \tau_1:Y_t \in A\}, T_A=\inf \{ t\geq 0:Y_t \notin A\},
\end{equation}
where $\tau_1=\inf \{ t\geq 0:Y_t \neq Y_0\}$.
\\
\indent The Green's function with respect to random walk $Z_t$ is defined by
\begin{equation}
\label{gr}    
g(x,y)=\int_0^\infty P_x^{\mathbb{Z}^d}[Z_t=y]dt, \text{ where } x,y\in \mathbb{Z}^d, \text{ and } g(\cdot)=g(0 , \cdot).
\end{equation}

\indent For any finite $K \subset \mathbb{Z}^d$ we define the equilibrium measure $e_K$ and the capacity $cap(K)$ by
\begin{equation}
\label{1.4}
e_K(x)=P_x^{\mathbb{Z}^d}[\tilde{H}_K=\infty]1_K(x) \text{ and } cap(K)=\sum_{x \in K} e_K(x).
\end{equation}
It is known that (see Proposition 6.5.2 and Proposition 6.5.6 in \cite{4}) there exist constants $c$ and $c'$ such that
\begin{equation}
\label{1.5}
cr^{d-2} \leq cap(B(0,r)) \leq c'r^{d-2} \text{ for } r\geq 1, d\geq 3.
\end{equation}
It is also known that (see Proposition 6.5.4, p138 in \cite{4}) for $d\geq 3$, there exist constants $c_1$, $c_2$ and $c_3$ such that for any $r\geq 1$, $K \subset B(0,r)$, 
\begin{equation}
\label{1.6}
c_1 \frac{e_K(y)}{cap(K)} \leq P_x^{\mathbb{Z}^d}[Z_{H_K}=y|H_K < \infty] \leq c_2 \frac{e_K(y)}{cap(K)} \text{ for }
y \in K, x\notin B(0,c_3 r). 
\end{equation}
\indent If $K \subset A \subset \mathbb{Z}^d$ with $K$ finite, we define the equilibrium measure and the capacity of $K$ relative to $A$ by
\begin{equation}
\label{1.7}
e_{K,A}(x)=P_x^{\mathbb{Z}^d}[\tilde{H}_K>T_A]1_K(x) \text{ and } cap_A(K)=\sum_{x \in K}e_{K,A}(x).
\end{equation}
Note that the following inequality follows directly from the definition of $e_K$ and $e_{K,A}$:
\begin{equation}
\label{1.12}
e_K(x) \leq e_{K,A}(x), \text{ for } x \in K \subset A \subset \mathbb{Z}^d.
\end{equation}
\indent  Similar to \cite{1}, we will need the following bounds on the probability of hitting time in $\mathbb{Z}^d$ and $\mathbb{T}_N$.
\begin{lemma}
Given $d \geq 3$ and moving mode matrix $Q$, there exists a constant $c>0$ such that 
\begin{equation}
\label{1.8}
P_x^{\mathbb{Z}^d}[H_{B(0,r_1)}<\infty] \leq c(r_1/r_2)^{d-2} \text{ for } 1\leq r_1 \leq r_2, x\notin B(0,r_2).
\end{equation}
\begin{equation}
\label{1.9}
\sup_{x \notin B(0,r_2)}P_x[H_{B(0,r_1)}<N^{2+\lambda}] \leq c(\lambda)(r_1/r_2)^{d-2} \text{ for } 1\leq r_1 \leq r_2\leq N^{1-3\lambda}, \lambda>0.
\end{equation}
\end{lemma}
\emph{Proof} (\ref{1.8}) follows from Proposition 6.5.1 in \cite{4} and (\ref{1.5}). The proof of (\ref{1.9}) is the same as the corresponding inequality (1.11) in \cite{1}. We omit the details here.\qed

Note that $Z_t$ is a martingale, so by the optional stopping theorem we have
\begin{lemma}
Given $d \geq 3$ and moving mode matrix $Q$, there exists a constant $c>0$ such that
\begin{equation}
\label{1.10}
P_x^{\mathbb{Z}^d}[T_{B(0,r_2)}<H_{B(0,r_1)}] \geq \frac{1}{r_2-r_1+c},\text{ for } r_1 < |x|_{\infty} \leq r_2.
\end{equation}
\end{lemma}
\indent Using the two lemmas above, we can prove the following bounds on equilibrium measures.
\begin{lemma}
Given $d \geq 3$ and moving mode matrix $Q$, there exists $c>0$ such that
\begin{equation}
\label{1.11}
e_K(x) \geq cr^{-1}, \text{ for } e_K(x)>0, \text{ where } K=B(0,r), r\geq 1.
\end{equation}
Furthermore  if  $r\geq 1, \lambda>0$, $K \subset B(0,r) \subset A = B(0,r^{1+\lambda})$, we have constants $c(\lambda)>0$ such that\begin{equation}
\label{1.13}
e_{K,A}(x) \leq (1+c(\lambda)r^{-\lambda})e_K(x), \text{ for all }x \in K.
\end{equation}
\end{lemma}
\emph{Proof} For a sufficiently large constant $c'$ we have $\inf_{x\notin B(0,c'r)}P_x^{\mathbb{Z}^d}[H_{B(0,r)}=\infty]\geq \frac{1}{2}$
(by(\ref{1.8})), so using (\ref{1.10}) in Lemma 1.2 we have (\ref{1.11}). For (\ref{1.13}) note that $e_{K,A}(x)=e_K(x)+P_x^{\mathbb{Z}^d}[T_A < \tilde{H}_K < \infty]$
and $P_x^{\mathbb{Z}^d}[T_A < \tilde{H}_K < \infty] \leq e_{K,A}(x)\sup_{y \notin A}P_y^{\mathbb{Z}^d}[H_K<\infty] \leq cr^{-\lambda}e_{K,A}(x)$(by (\ref{1.8})), so that $(1-cr^{-\lambda})e_{K,A}(x) \leq e_K(x)$ and (\ref{1.13}) is straightforward.\qed

We also need some notations for the trace of Poisson point processes of trajectories. Let $\Gamma=\Gamma(\mathbb{T}_N)$ or $\Gamma(\mathbb{Z}^d)$. When $\mu$ is a Poisson point process on $\Gamma$, we denote the trace of $\mu$ by
\begin{equation}
\label{1.14}
\mathcal{T}(\mu)=\cup_{\omega \in Supp(\mu)} \omega(0,\infty).
\end{equation}
If $\mu$ is a Poisson point process on $\Gamma \times [0,\infty)$, we denote the trace of $\mu$ up to label $u$ by
\begin{equation}
\label{1.15}
\mathcal{T}^u(\mu)=\mathcal{T}(\mu_u), \text{ where } \mu_u(dw)=\mu(dw \times [0,u]).
\end{equation}
If $\mu$ is a Poisson point process on the product space $\Gamma^n, n\geq 1$, we denote the trace of $\mu$ by
\begin{equation}
\mathcal{T}(\mu)=\cup_{(\omega_1, ..., \omega_n)\in Supp(\mu)} \cup_{j=1}^n \omega_j(0,\infty).
\end{equation}

The concept of random interlacements (which roughly speaking, are defined as a Poisson point process on the space of doubly infinite trajectories modulo time shift) was introduced by Sznitman in \cite{5}. Teixeira extended it to general transient weighted graphs in \cite{7} so that it still works in our setting. In fact, the famous \emph{random interlacements} $(\mathcal{I}^u)_{u \geq 0}$ is a family of random subsets of $\mathbb{Z}^d$ in a probability space $(\Omega_0,\mathcal{A}_0,Q_0)$. In this article we need the following facts from \cite{5} and \cite{7} (see e.g. (2.11) and (2.23) in \cite{7}).
\begin{proposition}
\label{ri}
The following statements hold for random interlacements $(\mathcal{I}^u)_{u \geq 0}$.\\
(i). $(\mathcal{I}^u \cap K)_{u \geq 0} \overset{\text{law}}{=}(\mathcal{T}^u(\mu_K) \cap K)_{u \geq 0}$ for all finite $K \subset \mathbb{Z}^d$ where $\mu_K$ is a Poisson point process on $\Gamma(\mathbb{Z}^d) \times [0,\infty)$ of intensity $P_{e_K}^{\mathbb{Z}^d}  \otimes \lambda$, here $\lambda$ denotes the Lebesgue measure.\\
(ii). The law of $\mathcal{I}^u$ under $Q_0$ is translation-invariant for all $u \geq 0$.\\
(iii). $\mathcal{I}^u$ is increasing with respect to $u$ in the sense that $Q_0-$almost surely $\mathcal{I}^v \subset \mathcal{I}^u$ for $v \leq u$. \\
(iv). If $\mathcal{I}_1^u$ and $\mathcal{I}_2^v$ are independent with the laws of $\mathcal{I}^u$ and $\mathcal{I}^v$ under $Q_0$ respectively, then $(\mathcal{I}_1^u, \mathcal{I}_1^u \cup \mathcal{I}_2^v)$ has the law of $(\mathcal{I}^u,\mathcal{I}^{u+v})$  under $Q_0$.
\end{proposition}
\indent For the law of random interlacements we also have (see Theorem 4.1.1, p75 in \cite{4})
\begin{equation}
\label{1.18}
Q_0[x\notin \mathcal{I}^u] = \exp(- \frac{u}{g(0)}) \text{ and } Q_0[x,y \notin \mathcal{I}^u] = \exp(- \frac{2u}{g(0)+g(x-y)}).
\end{equation}
\indent Lemma 1.5 in \cite{8} gives a strong estimation for the sum of the second term in (\ref{1.18}), and can actually be straightforwardly generalized into finite range, symmetric, irreducible random walks. In fact we have
\begin{lemma}
Given $d \geq 3$ and moving mode matrix $Q$, there exists a $c_4>1$ and some $c>0$ such that for any $K \subset \mathbb{Z}^d$ with $0 \notin K$ and $u \geq 0$
\begin{equation}
\label{1.19}
\sum_{v\in K} Q_0[0, v \notin \mathcal{I}^u] \geq c|K| (Q_0[0 \notin \mathcal{I}^u])^2(1+u)+ce^{-c_4\frac{u}{g(0)}}.
\end{equation}
\end{lemma}

Finally we recall that the cover time $T_C^F$ of a set $F \subset \mathbb{T}_N$ is defined by
\begin{equation}
\label{1.20}
T_C^F:= \max_{x\in F}H_x,
\end{equation}
and sometimes we use $T_C$ to represent $T_C^{\mathbb{T}_N}$. 

As in \cite{1}, we'll use the notation
\begin{equation}
\label{1.21}
u_F(z)=g(0)(\log|F|+z),
\end{equation}
so that the event $\{\frac{T_C^F}{g(0)N^d}- \log |F| \leq z\} = \{ T_C^F \leq u_F(z) N^d\}$.

\section{From random walk to random  interlacements}

In this section we state a powerful coupling result between the random walk and random interlacements  which is crucial in the proof of Theorem \ref{T1} and then present the proof of it. This result is also a generalization of Theorem 2.2 in \cite{1}.
\subsection{Statement of the coupling result}
\indent For $n \geq 1$, and $x_1, ..., x_n \in \mathbb{T}_N$ we define the separation $s=s(x_1, ..., x_n)$ of the vertices $x_1, ..., x_n$ by 
\begin{equation}
\label{2.23}
s=s(x_1, ..., x_n)=\left\{
\begin{aligned}
&N & \text{ if } n=1,\\
&\min_{i \neq j}d_{\infty}(x_i,x_j) & \text{ if } n>1.
\end{aligned}
\right.
\end{equation}
\indent From now on, we need an arbitrarily small constant $\epsilon_0>0$ which does not depend on $N$, and we define the box 
\begin{equation}
\label{2.24}
\Xi^A=B(0,s^{1-\epsilon_0}).
\end{equation}
Note that for convenience we sometimes view $\Xi^A$ (or such boxes) also as $\{L(x):x\in \Xi^A\}$, a subset of $\mathbb{T}_N$, although it is in fact a subset of $\mathbb{Z}^d$ from the definition.

The following result will couple the trace of random walk in the boxes $\Xi^A+x_1, ..., \Xi^A+x_n$ with independent random interlacements. To be precise we have
\begin{proposition}[a generalization of Theorem 2.2 in \cite{1}]
\label{T2}
\indent Given $d \geq 3$ and moving mode matrix $Q$, let $n \geq 1$, $x_1, ..., x_n \in \mathbb{T}_N$ be distinct from each other and have separation $s$ (recall (\ref{2.23}) for the definition of $s$), $\epsilon_0 \in (0,1)$, then there exist constants $c>0$ and $c_5=c_5(\epsilon_0)>0$ small enough, such that if $u \geq s^{-c_5}$, $1 \geq \delta \geq \frac{1}{c_5} s^{-c_5}$, $n \leq s^{c_5}$, we can construct a space $(\Omega_1,\mathcal{A}_1,Q_1)$ with a random walk $Y_{\cdot}$ with law $P$ and $n$ independent copies of random interlacements $(\mathcal{I}_i^v)_{v \geq 0}, i= 1, ..., n$, each with the law of $(\mathcal{I}^v)_{v \geq 0}$ under $Q_0$, and
\begin{equation}
\label{2.25}
Q_1[\mathcal{I}_i^{u(1-\delta)} \cap \Xi^A \subset (Y(0, uN^d)- x_i) \cap \Xi^A \subset \mathcal{I}_i^{u(1+\delta)} \cap \Xi^A] \geq 1-cue^{-c s^{c_5}} \text{ for all } i.
\end{equation}
\end{proposition}

The following corollary could be derived 	directly from Proposition \ref{T2}.
\begin{corollary}
Given $d \geq 3$ and moving mode matrix $Q$, there exist constants $c, c', c_6>0$  such that if $N^{-c_6} \leq u \leq N^{c_6}$, then for all $x \in \mathbb{T}_N$
\begin{equation}
\label{2.26}
Q_0[0 \notin \mathcal{I}^u](1-c'N^{-c}) \leq P(x \notin Y(0,uN^d)) \leq Q_0[0 \notin \mathcal{I}^u](1+c' N^{-c}).
\end{equation}
\end{corollary}
\emph{Proof} Taking $n=1, x_1=x, \delta = \frac{1}{c_5} N^{-c_5}$ in Proposition \ref{T2} and using (\ref{1.18}), we get the inequality in (\ref{2.26}) straightforwardly.\qed

\subsection{The coupling lemmas}

We now concentrate on the proof of Proposition \ref{T2} where the proof follows an argument similar to that in \cite{1}. More precisely, we will first introduce several coupling lemmas and show that Proposition \ref{T2} could be derived from them. The proof of these lemmas are postponed in the subsequent subsections.

For the rest of the paper we assume that centers of boxes are given: 
\begin{equation}
\label{4.30}
x_1,..,x_n \in \mathbb{T}_N \text{ whose separation is } s.
\end{equation}
We also define the concentric boxed $\Xi^B \subset \Xi^C$ around $\Xi^A$ (recall (\ref{2.24}) for the definition of $\Xi^A$) by
\begin{equation}
\label{4.31}
\Xi^A\subset \Xi^B=B(0,s^{1-\frac{\epsilon_0}{2}}) \subset \Xi^C=B(0,s^{1-\frac{\epsilon_0}{4}}).
\end{equation}
For convenience we introduce the notation
\begin{equation}
\label{4.32}
	\bar{M}=\cup_{i=1}^{n}M_i \text{ where } M_i=M+x_i \text{ for any } M\subset \mathbb{T}_N.
\end{equation} 
Note that for any $0<\epsilon_0<1$ there exists some $c(\epsilon_0)$ such that if $s \geq c(\epsilon_0)$ then all the $\Xi^C_i$ are disjoint. \\
\indent We also introduce $U$, the first time random walk spends a long time outside of $\overline{\Xi^C}$ (roughly speaking a time long enough to mix), defined by
\begin{equation}
\label{4.33}
U=\inf\{t \geq t^{*}:Y(t-t^{*},t)\cap \overline{\Xi^C}=\emptyset\} \text{ where } t^{*}=N^{2+\frac{\epsilon_0}{100}}.
\end{equation}
And we introduce the intensity measures $\kappa_1$, $\kappa_2$ on $\Gamma(\mathbb{T}_N)$ and $\kappa_3$ on $\Gamma(\mathbb{Z}^d)$ defined by
\begin{equation}
\label{4.34}
\kappa_1(W)=P_{e}[Y_{\cdot \wedge U} \in W] \text{ where } e(x)=\sum_{i=1}^n e_{\Xi^A_i}(x).
\end{equation}
\begin{equation}
\label{4.35}
\kappa_2(W)=P_{e}[Y_{\cdot \wedge T_{\overline{\Xi^B}}} \in W].
\end{equation}
\begin{equation}
\label{4.36}
\kappa_3(W)=P_{e_{\Xi^A}}^{\mathbb{Z}^d}[Z_{\cdot \wedge T_{\Xi^B}} \in W].
\end{equation}
\indent We now state the three coupling lemmas.
\begin{lemma}[a generalization of  Lemma 4.1 in \cite{1}]
\label{L3.1}
\indent Given $d \geq 3$ and moving mode matrix $Q$, there exist some $c, c(\epsilon_0)>0$ such that if $s \geq c(\epsilon_0)$, $u \geq s^{-c(\epsilon_0)}$, $1 \geq \delta \geq cs^{-c(\epsilon_0)}$, and $n \leq s^{c(\epsilon_0)}$ we can construct a probability space $(\Omega, \mathcal{A}, Q)$ with a random walk $Y_{t}$ with law $P$ and independent Poisson point processes $\mu_1,\mu_2$, on $\Gamma(\mathbb{T}_N)$ such that $\mu_1$ has intensity $u(1-\delta)\kappa_1$, $\mu_2$ has intensity $2u\delta \kappa_1$, and $Q[I_1] \geq 1-cue^{-cs^{c(\epsilon_0)}}$, where
\begin{equation}
\label{4.37}
I_1=\{\mathcal{T}(\mu_1) \cap \overline{\Xi^A} \subset Y(0,uN^d) \cap \overline{\Xi^A} \subset \mathcal{T}(\mu_1+\mu_2) \cap \overline{\Xi^A}\}.
\end{equation}
\end{lemma}
\begin{lemma}[a generalization of  Lemma 4.2 in \cite{1}]
\label{L3.2}
\indent Given $d \geq 3$ and moving mode matrix $Q$, there exist some $c, c(\epsilon_0)>0$ such that if $s \geq c(\epsilon_0)$, $u \geq s^{-c(\epsilon_0)}$, $1 \geq \delta \geq cs^{-c(\epsilon_0)}$, $n \leq s^{c(\epsilon_0)}$ then we can construct a probability space $(\Omega,\mathcal{A},Q)$ with Poisson point processes $\nu$, $\nu_1$, $\nu_2$ on $\Gamma(\mathbb{T}_N)$ such that $\nu$ has intensity $u\kappa_1$, $\nu_1$ has intensity $u(1-\delta)\kappa_2$, $\nu_2$ has intensity $2u\delta \kappa_2$, $\nu_1$ and $\nu_2$ are independent to each other, $Q[I_2]=1$ and $Q[I_3] \geq  1-ce^{-cs^{c(\epsilon_0)}}$, where
\begin{equation}
\label{4.38}
I_2=\{\mathcal{T}(\nu_1) \cap \overline{\Xi^A} \subset \mathcal{T}(\nu) \cap \overline{\Xi^A}\},
\end{equation}
\begin{equation}
\label{4.39}
I_3=\{\mathcal{T}(\nu) \cap \overline{\Xi^A} \subset \mathcal{T}(\nu_1+\nu_2) \cap \overline{\Xi^A}\}.
\end{equation}
\end{lemma}
\begin{lemma}[a generalization of  Lemma 4.3 in \cite{1}]
\label{L3.3}
\indent Given $d \geq 3$ and moving mode matrix $Q$, there exist some $c, c(\epsilon_0)>0$ such that if $s \geq c(\epsilon_0)$, $u \geq 0$, $1 \geq \delta \geq cs^{-c(\epsilon_0)}$ then we can construct a probability space $(\Omega,\mathcal{A},Q)$ with a Poisson point process $\eta$  on $\Gamma(\mathbb{Z}^d)$ with intensity measure $u\kappa_3$ and independent random sets $\mathcal{I}_1,\mathcal{I}_2\subset \mathbb{Z}^d$ such that $\mathcal{I}_1$ has the law of $\mathcal{I}^{u(1-\delta)}$, $\mathcal{I}_2$ has the law of $\mathcal{I}^{2u\delta}$ and $Q[I_4] \geq  1-ce^{-cs^{c(\epsilon_0)}}$ (see Proposition \ref{ri} for the definition of random interlacements), where
\begin{equation}
\label{4.40}
I_4=\{\mathcal{I}_1 \cap \Xi^A \subset \mathcal{T}(\eta) \cap \Xi^A \subset (\mathcal{I}_1 \cup \mathcal{I}_2) \cap \Xi^A\}.
\end{equation}
\end{lemma}
\emph{Proof of Proposition \ref{T2}} Take $\nu=\mu_1$ and $\nu=\mu_1+\mu_2$ respectively, and note that $\kappa_2$ start from each $\Xi^A_i$ is equivalent to $\kappa_3$. We could combine the three couplings in the lemmas above to get Proposition \ref{T2} straightforwardly.\qed\\

\indent Proposition \ref{T2} has now been reduced to Lemma \ref{L3.1}, Lemma \ref{L3.2} and Lemma \ref{L3.3}.\\
\indent Lemma \ref{L3.1} will be proved in Section 3.3 and Lemma \ref{L3.2} will be proved in Section 3.4. The proof of Lemma \ref{L3.3} is the same as the proof of Proposition 4.4 in \cite{1}. We omit it here.
\subsection{Quasistationary Distribution}
\indent In this section we first introduce the \emph{quasistationary distribution}, and use it to prove Lemma \ref{L3.1}. \\
\indent We define the $(N^d-|\overline{\Xi^C}|) \times (N^d -|\overline{\Xi^C}|)$ matrix $(P^{\overline{\Xi^C}})_{x,y\in{\mathbb{T}_N \setminus \overline{\Xi^C}}}=P_x[Y_{\tau_1}=y]1_{x\neq y}$ to denote the probability transition matrix for the random walk restricted in $\mathbb{T}_N \setminus \overline{\Xi^C}$. We point out that given $d\geq 3$ and the moving mode matrix $Q$, $\mathbb{T}_N \setminus \overline{\Xi^C}$ is connected by the random walk restricted in $\mathbb{T}_N \setminus \overline{\Xi^C}$ when $s$ is greater than some $c(\epsilon_0)>0$. 

To show this we first note that given $d\geq 3$ and the moving mode matrix $Q$, there exists $s_0>0$ such that $0$ is connected to $\{e_1, e_2, ..., e_d\}$ (see (\ref{ei}) for the definition) by the random walk restricted in $B(0, s_0)$. So that $\{x\in \mathbb{T}_N: |x- \overline{\Xi^C}|_{\infty} >s_0\}$ is connected by the random walk restricted in $\mathbb{T}_N \setminus \overline{\Xi^C}$. As a result, we could simply choose one suitable step from the moving mode matrix $Q$ and repeat this step to further show that each point in $\mathbb{T}_N \setminus \overline{\Xi^C}$ is connected to $\{x\in \mathbb{T}_N: |x- \overline{\Xi^C}|_{\infty} >s_0\}$ by the random walk restricted in $\mathbb{T}_N \setminus \overline{\Xi^C}$ when $s$ is greater than some $c(\epsilon_0)>0$.

Perron-Frobenius theorem (Theorem 8.2 in \cite{9}) implies that the real symmetric, non-negative and irreducible matrix $P^{\overline{\Xi^C}}$ has a unique largest eigenvalue $\lambda_1^{\overline{\Xi^C}}$ with a non-negative normalized eigenvector $v_1$. We let $\lambda_2^{\overline{\Xi^C}}$ denote the second largest eigenvalue of $P^{\overline{\Xi^C}}$. The quasistationary distribution $\sigma$ on $\mathbb{T}_N \setminus \overline{\Xi^C}$ is then defined by
\begin{equation}
\label{5.41}
\sigma(x)=\frac{(v_1)_x}{v_1^{T}\boldsymbol{1}} \text{ for } x\in \mathbb{T}_N \setminus \overline{\Xi^C}.
\end{equation}
\indent Since $\mathbb{T}_N\setminus \overline{\Xi^C}$ is connected (when $s \geq c(\epsilon_0)$) it holds that (see (6.6.3) in \cite{10})
\begin{equation}
\label{5.42}
\lim_{t \to \infty} P_x[Y_t=y|H_{\overline{\Xi^C}}>t] =\sigma(y) \text{ for all } x,y \in \mathbb{T}_N \setminus \overline{\Xi^C}.
\end{equation}
\indent We first need the following lemma on the lower bound of $\sigma(x)$, the proof is slightly different from the corresponding lemma in \cite{1}.
\begin{lemma}[a generalization of  Lemma 5.4 in \cite{1}]
\label{L4.1}
\indent Given $d \geq 3$ and moving mode matrix $Q$, there exist $c, c(\epsilon_0)>0$ such that if $ s\geq c(\epsilon_0)$ we have
\begin{equation}
\label{5.48}
\inf_{x\in \mathbb{T}_N \setminus \overline{\Xi^C}} \sigma(x) \geq N^{-cn}.
\end{equation}
\end{lemma}
\indent \emph{Proof } We sometimes directly use $c$ to represent different positive constants which means that there exists such $c>0$. We first define concentric boxes $\Xi_D,\Xi_E$ and $\Xi_F$ by
\begin{equation}
\label{5.49}
\Xi_D=B(0,s^{1-\frac{\epsilon_0}{8}}) \subset \Xi_E=B(0,s^{1- \frac{\epsilon_0}{16}}) \subset \Xi_F=B(0,s^{1-\frac{\epsilon_0}{32}}).
\end{equation}
\indent Let $y$ be the maximum of $\sigma(\cdot)$. Since $\sigma(\cdot)$ is a probability distribution we have $\sigma(y) \geq N^{-d}$.
Also by reversibility we have for any $x \notin \overline{\Xi^C}$ and $t\geq 0$ that $P_x[Y_t=y, H_{\overline{\Xi^C}}>t]=P_y[Y_t=x, H_{\overline{\Xi^C}}>t]$ and thus we have
\begin{equation}
\label{5.50}
P_x[Y_t=y|H_{\overline{\Xi^C}}>t]=P_y[Y_t=x|H_{\overline{\Xi^C}}>t]\frac{P_x[H_{\overline{\Xi^C}}>t]}{P_y[H_{\overline{\Xi^C}}>t]}.
\end{equation}
\indent Since by the strong Markov property $P_x[H_{\overline{\Xi^C}}>t] \geq P_x[H_y < H_{\overline{\Xi^C}}]P_y[H_{\overline{\Xi^C}} >t]$ we see, by taking the limit $t \to \infty$ in (\ref{5.50}) and using (\ref{5.42}), that $\sigma(y) \geq \sigma(x)P_x[H_y < H_{\overline{\Xi^C}}] \geq N^{-d}P_x[H_y < H_{\overline{\Xi^C}}]$. To prove (\ref{5.48}) it thus suffices to show that
\begin{equation}
\label{5.51}
P_x[H_y < H_{\overline{\Xi^C}}] \geq N^{-cn} \text{ for all } x,y \notin \overline{\Xi^C}.
\end{equation}
\indent If $x \in \Xi^D_i\setminus \Xi^C_i$ for some $i=1,2, ..., n$, it follows from (\ref{1.10}) that $P_x[T_{\Xi^D_i} < H_{\overline{\Xi^C}}] \geq c N^{-1}$, so that by the Markov property $P_x[H_y < H_{\overline{\Xi^C}}] \geq c N^{-1}\inf_{x' \notin \Xi^D_i} P_{x'}[H_y < H_{\overline{\Xi^C}}]$. If $x \notin \overline{\Xi^D}$ and $y \in \overline{\Xi^F}\setminus \overline{\Xi^C}$ then we can use that by reversibility $P_x[H_y <H_{\overline{\Xi^C}}]=P_y[H_x <H_{\overline{\Xi^C}}] $ and (\ref{1.10}) to show that $P_x[H_y<H_{\overline{\Xi^C}}] \geq c N^{-1}\inf_{y' \notin \overline{\Xi^F}}P_x[H_{y'} <H_{\overline{\Xi^C}}]$. So we only need to show
\begin{equation}
\label{5.52}
P_x[H_y<H_{\overline{\Xi^C}}] \geq N^{-cn} \text{ for all } x \notin \overline{\Xi^D}, y \notin \overline{\Xi^F}.
\end{equation}
\indent Now fix $y \notin \overline{\Xi^F}$ and note that $P_x[H_y < H_{\overline{\Xi^C}}] \geq cs^{-c}$ for all $x \in \partial_e(\Xi^D+y) \cup (\Xi^D+y)$, since $y+\Xi^E \cap \overline{\Xi^C}=\emptyset$ when $s$ is greater than some $c(\epsilon_0)>0$ and $P_x[H_y < T_{\Xi^E+y}] \geq cs^{-c}$(by Proposition 6.3.2 and 6.3.5 in \cite{4}), where $\partial_e(V):=\{x \notin V: |x-V|_2=1\}$. Therefore to prove (\ref{5.52}) it suffices to show
\begin{equation}
\label{5.53}
P_{x_1}[H_y<H_{\overline{\Xi^C}}] \geq N^{-cn}P_{x_2}[H_y <H_{\overline{\Xi^C}}] \text{ for all } x_1, x_2 \notin \overline{\Xi^D} \cup (\Xi^D+y).
\end{equation}
\indent Consider the function $x \to P_x[H_y<H_{\overline{\Xi^C}}]$. The function is non-negative and harmonic w.r.t the random walk $Z_t$ in $(\overline{\Xi^C}\cup \{y\})^{c}$. Thus by the Harnack inequality (Theorem 6.3.9 in \cite{4}) we have that, for any $z \in \mathbb{T}_N$ and $r\geq 0$ for which $B(z,2(r+1)) \cap (\overline{\Xi^C} \cup \{y\}) = \emptyset$ (without loss of generality, we could assume that $B(z,2(r+1))$ doesn't touch the boundary of torus),
\begin{equation}
\label{5.54}
\inf_{x \in B(z,r+1)}P_x[H_y<H_{\overline{\Xi^C}}] \geq c \sup_{x \in B(z,r+1)}P_x[H_y<H_{\overline{\Xi^C}}].
\end{equation}
\indent Lemma 5.5 in \cite{1} shows that if $s \geq c(\epsilon_0)$, one can cover $(\overline{\Xi^D}\cup (\Xi^D+y))^c$ by $m\leq cn \log N$ balls $B(z_i, r_i), i = 1, ...,n$ that satisfy $B(z_i, 2(r_i + 1)) \cap (\overline{\Xi^C} \cup {y})= \emptyset$, so that by iterating this inequality we get the proof of Lemma \ref{L4.1}.\qed \\
\indent We now give an estimate on the rate of convergence in (\ref{5.42}).
\begin{lemma}[a generalization of  Proposition 5.1 in \cite{1}]
\label{L4.2}
\indent Given $d \geq 3$ and moving mode matrix $Q$, there exist some $c, c(\epsilon_0)>0$ such that if $n \leq s^{c(\epsilon_0)}$ and $ s\geq c(\epsilon_0)$ then
\begin{equation}
\label{5.43}
\sup_{x,y \in \mathbb{T}_N \setminus\overline{\Xi^C}}|P_x[Y_{t^*}=y|H_{\overline{\Xi^C}}>t^*]-\sigma(y)| \leq ce^{-cN^{c(\epsilon_0)}}.
\end{equation}
\end{lemma}
\emph{Proof } Note that for $V\in \overline{\Xi^C}$ we could generalize Proposition 3.2 in \cite{12} for finite range, symmetric, irreducible random walks in $\mathbb{T}_N$ to get:
\begin{equation}
\label{DV}
\mathcal{D}(g^*,g^*)(1-2\sup_{x\notin \overline{\Xi^D}}|f^*(x)|)\leq \frac{1}{\mathbb{E}H_V}\leq\mathcal{D}(g^*,g^*)\pi(\mathbb{T}_N \setminus\overline{\Xi^D})^{-2}
\end{equation}
where $\pi$ is the invariant distribution (or uniform distribution in our setting), $\mathcal{D}$ is the Dirichlet form given by
\begin{equation}
\mathcal{D}(f,g):=-\frac{1}{2}\sum_{x,y\in \mathbb{T}_N}(f(x)-f(y))(g(x)-g(y))\pi_x \frac{q_{xy}}{q_{xx}}
\end{equation}
for any real-valued function $f,g$ on $\mathbb{T}_N$,
\begin{equation}
f^*(x)=1-\frac{\mathbb{E}_x H_V}{\mathbb{E} H_V}
\end{equation}
and
\begin{equation}
g^*(x)= \mathbb{P}_x(H_V\leq H_{\mathbb{T}_N \setminus\overline{\Xi^D}}).
\end{equation}
To prove this inequality we need to take $A=V$, $C=\mathbb{T}_N \setminus \overline{\Xi^D}$ in Proposition 3.2 in \cite{12}. For the LHS of (\ref{DV}) we simply replace all the $\partial_i C$ by $C$ in the proof and for the RHS, since the original proof of Lemma 3.1 in \cite{12} is actually done on general weighted graphs, it naturally applies to our problem. Thus Proposition 3.3 in \cite{22} and Lemma 5.2 in \cite{1} could in turn be generalized and we'll have:

Given $d \geq 3$ and moving mode matrix $Q$, there exist some $c(\epsilon_0)>0$ such that if $s\geq c(\epsilon_0)$, $V \subset \overline{\Xi^C}$, then let $V_i =(V \cap \Xi^C_i)-x_i \subset \mathbb{Z}^d, i=1, ..., n$, we have
\begin{equation}
\label{5.44}
\frac{N^d}{E[H_V]\sum_{i=1}^{n}cap(V_i)} \leq 1+c(\epsilon_0)s^{-c(\epsilon_0)}.
\end{equation}
Also there exists some $c(\epsilon_0)>0$ such that if $V \subset \overline{\Xi^B}$ and $n \leq s^{c(\epsilon_0)}$ then
\begin{equation}
\label{5.45}
1-c(\epsilon_0)s^{-c(\epsilon_0)}\leq \frac{N^d}{E[H_V]\sum_{i=1}^{n}cap(V_i)},
\end{equation}
and
\begin{equation}
\label{5.46}
(1-c(\epsilon_0)s^{-c(\epsilon_0)})E[H_V] \leq \inf_{x \notin \overline{\Xi^C}}E_x[H_V] \leq \sup_{x\in \mathbb{T}_N}E_x[H_V] \leq (1+c(\epsilon_0)s^{-c(\epsilon_0)})E[H_V].
\end{equation}
\indent With this result we could generalize Lemma 5.3 in \cite{1} to have that, given $d \geq 3$ and moving mode matrix $Q$, there exist some $c, c(\epsilon_0)>0$ such that if $n \leq s^{c(\epsilon_0)}$ and $s \geq c(\epsilon_0)$ we have
\begin{equation}
\label{5.47}
\lambda_1^{\overline{\Xi^C}}-\lambda_2^{\overline{\Xi^C}} \geq cN^{-2}.
\end{equation}
Thus the remaining proof of Lemma \ref{L4.2} follows the same arguments as in the proof of Proposition 5.1 in \cite{1}. And we omit the details.\qed \\
\indent Finally, we give the proof of Lemma \ref{L3.1} as follows.\\\\
\indent \emph{Proof of Lemma \ref{L3.1} } As Lemma 5.6 in \cite{1} we have, given $d \geq 3$ and moving mode matrix $Q$, there exist some $c, c(\epsilon_0)>0$ such that if $s \geq c(\epsilon_0)$ and $n \leq s^{c(\epsilon_0)}$ then for all $i=1, ..., n$,
\begin{equation}
\label{5.55}
\frac{e_{\Xi^A}(x-x_i)}{n\text{cap}(\Xi^A)}(1-cs^{-c(\epsilon_0)}) \leq P_{\sigma}[Y_{H_{\overline{\Xi^A}}}=x] \leq \frac{e_{\Xi^A}(x-x_i)}{n\text{cap}(\Xi^A)}(1+cs^{-c(\epsilon_0)}) \text{ for all } x \in \Xi^A_i.
\end{equation}
So that with Lemma \ref{L4.2}, we could straightforwardly generalize Proposition 4.1 in \cite{1} into Lemma \ref{L3.1}.\qed

\subsection{From the Torus to $\mathbb{Z}^d$} 
\indent In this section we will prove Lemma \ref{L3.2}.\\
\indent First we define for $\omega \in \Gamma(\mathbb{T}_N)$ the successive returns $\hat{R}_k=\hat{R}_k(\omega)$ to $\overline{\Xi^A}$ and departures $\hat{D}_k=\hat{D}_k(\omega)$ from $\overline{\Xi^B}$ as follows
\begin{equation}
\label{6.56}
\hat{R}_1=H_{\overline{\Xi^A}}, \hat{R}_k=\inf\{ t \geq \hat{D}_{k-1}: \omega_t \in \overline{\Xi^A}\}, \hat{D}_k=\inf\{t \geq \hat{R}_k: \omega_t \notin \overline{\Xi^B} \}, k\geq 1.
\end{equation}
\indent To extract the successive visits to $\overline{\Xi^A}$ of an excursion we furthermore define for each $i \geq 1$ the map $\phi_i$ from $\{\hat{R}_i<U<\hat{R}_{i+1}\} \subset \Gamma(\mathbb{T}_N)$ into $\Gamma(\mathbb{T}_N)^{i}$ by
\begin{equation}
\label{6.57}
(\phi_i(\omega))_j=\omega((\hat{R}_j + \cdot) \wedge \hat{D}_j) \text{ for } j=1,..,i, \omega \in \{\hat{R}_i<U<\hat{R}_{i+1}\} \subset \Gamma(\mathbb{T}_N), i\geq 1.
\end{equation}
\indent For each $i \geq 1$ we will apply this map to the Poisson point process $1_{\hat{R}_i<U<\hat{R}_{i+1}}\nu$ to get Poisson point process $\mu_i$ of intensity $u\kappa^{i}_{1}$ (recall Lemma \ref{L3.2} for the definition of $\nu$) on $\Gamma(\mathbb{T}_N)^{i}$, where
\begin{equation}
\label{6.58}
\kappa^{i}_1 = \phi_i \circ (1_{\hat{R}_i<U<\hat{R}_{i+1}}\kappa_1), i\geq 1.
\end{equation}
\indent To prove Lemma \ref{L3.2}, we first need the following bound on $\kappa_1$.
\begin{lemma}[a generalization of Lemma 7.1 in \cite{1}]
\label{L5.1}
\indent Given $d \geq 3$ and moving mode matrix $Q$, there exist some $c, c(\epsilon_0) >0$ such that if $s \geq c(\epsilon_0)$ and $n \leq s^{c(\epsilon_0)}$ we have for any measurable $W\subset \Gamma(\mathbb{T}_N)$,
\begin{equation}
\label{6.60}
(1-cs^{-c(\epsilon_0)})\kappa_2(W) \leq \kappa_1^1(W) \leq \kappa_2(W).
\end{equation}
\end{lemma}
\indent \emph{Proof } Recalling that $\kappa_1^1(W)=P_e[Y_{\cdot \wedge T_{\overline{\Xi^B}}} \in W, U< \hat{R}_2]$, so the upper bound follows directly from the definition of $\kappa_2$ (see \ref{4.35} for the definition of $\kappa_2$). Furthermore $\kappa_1^1(W) \geq \kappa_2(W)\inf_{x \notin \overline{\Xi^B}}P_x[H_{\overline{\Xi^A}}>T_{\overline{\Xi^D}}]\inf_{x \notin \overline{\Xi^D}}P_x[H_{\overline{\Xi^A}}>U]$ by the strong Markov property. Note that there exist some $c, c(\epsilon_0)>0$ such that when $s \geq c(\epsilon_0)$ we have $\inf_{x \notin \overline{\Xi^B}}P_x[H_{\overline{\Xi^A}} >T_{\overline{\Xi^D}}] = \inf_{x \notin \Xi^B}P_x^{\mathbb{Z}^d}[H_{\Xi^A}>T_{\Xi^D}]\geq 1-cs^{-c(\epsilon_0)}$, where the last inequality comes from (\ref{1.8}). On the other hand there exist some $c(\epsilon_0)>0$ such that when $n \leq s^{c(\epsilon_0)}$, we have
\begin{equation}
\label{6.61}
\sup_{x \notin \overline{\Xi^D}}P_x[H_{\overline{\Xi^C}}<U]=\sup_{x \notin \overline{\Xi^D}}P_x[H_{\overline{\Xi^C}}<N^{2+\frac{\epsilon_0}{100}}]\leq c(\epsilon_0)s^{-c(\epsilon_0)}.
\end{equation}
where the last inequality comes from (\ref{1.9}), so that $\inf_{x \notin \overline{\Xi^D}}P_x[H_{\overline{\Xi^A}}>U]\geq 1-c(\epsilon_0)s^{-c(\epsilon_0)}$ and the lower bound follows.\qed \\
\indent The following lemma is also crucial in the proof of Lemma \ref{L3.2}.
\begin{lemma}[a generalization of Lemma 7.5 in \cite{1}]
\label{L5.2}
\indent Given $d \geq 3$ and moving mode matrix $Q$, there exist some $c, c(\epsilon_0)>0$ such that if $ s\geq c(\epsilon_0)$ and $n \leq s^{c(\epsilon_0)}$ we have
\begin{equation}
\label{6.63}
h_x(y):=P_x[H_{\overline{\Xi^A}}<U, Y_{H_{\overline{\Xi^A}}}=y] \leq s^{-\frac{\epsilon_0}{4}}\bar{e}(y)
\end{equation}
for all $x\notin \overline{\Xi^B}$, $y\in \overline{\Xi^A}$, where  $\bar{e}=\frac{e}{ncap(\Xi^A)}$ denotes the normalization of the measure $e$.
\end{lemma}
\emph{Proof } By the strong Markov property we can assume that $x \in \overline{\Xi^C}$. We have
\begin{equation}
\label{6.64}
h_x(y)\leq P_x[H_{\overline{\Xi^A}}<T_{\overline{\Xi^D}}, Y_{H_{\overline{\Xi^A}}}=y]+P_x[T_{\overline{\Xi^D}}<H_{\overline{\Xi^A}}<U,Y_{H_{\overline{\Xi^A}}}=y].
\end{equation}
By the strong Markov property and (\ref{6.61}) we have that there exist some $c(\epsilon_0)>0$ such that for $s\geq c(\epsilon_0)$,
\begin{equation}
\label{6.65}
\begin{aligned}
P_x[T_{\overline{\Xi^D}}<H_{\overline{\Xi^A}}<U,Y_{H_{\overline{\Xi^A}}}=y] & \leq \sup_{z\notin \overline{\Xi^D}}P_z[H_{\overline{\Xi^C}}<U]\sup_{z\in \overline{\Xi^C}\setminus \overline{\Xi^B}}h_z(y) \\
& \leq c(\epsilon_0) s^{-c(\epsilon_0)}\sup_{z\in \overline{\Xi^C}\setminus \overline{\Xi^B}}h_z(y).
\end{aligned}
\end{equation}
So that if we take $x$ in (\ref{6.64}) to be the maximal point of $h_x(y)$ on $\overline{\Xi^C}\setminus \overline{\Xi^B}$, we have that there exist some $c(\epsilon_0)>0$ such that for $s \geq c(\epsilon_0)$, 
\begin{equation}
\label{6.66}
\begin{aligned}
\sup_{z\in \overline{\Xi^C}\setminus \overline{\Xi^B}}h_z(y) &\leq c(\epsilon_0) \sup_{z\in \overline{\Xi^C}\setminus\overline{\Xi^B}}P_z[H_{\overline{\Xi^A}}<T_{\overline{\Xi^D}}, Y_{H_{\overline{\Xi^A}}}=y]\\
&\leq c(\epsilon_0)\sup_{z\in {\Xi^C_i}\setminus {\Xi^B_i}}P_z^{\mathbb{Z}^d}[H_{{\Xi^A_i}}<\infty, Z_{H_{{\Xi^A_i}}}=y].
\end{aligned}
\end{equation}
Here we suppose that $y\in \Xi^A_i \in \overline{\Xi^A}$. Using (\ref{1.6}) and (\ref{1.8}) in (\ref{6.66}), we can get the inequailty in (\ref{6.63}) (Note that there might be a mistake in the proof of Lemma 7.5 in \cite{1} since $h_x(y)$ seems not to be a harmonic function).\qed \\\\
\indent At the end of this section we conclude the proof of Lemma \ref{L3.2}.\\
\indent \emph{Proof of Lemma \ref{L3.2} } By Lemma \ref{L5.2} above we could generalize Lemma 7.3 in \cite{1} straightforwardly to get that there exist some $c, c(\epsilon_0)>0$ such that if $s \geq c(\epsilon_0)$ and $n \leq s^{c(\epsilon_0)}$ then for all $i \geq 2$
\begin{equation}
\label{6.62}
\begin{aligned}
& \kappa_1^i \leq \tilde{\kappa}_1^{i}, \text{ where } \tilde{\kappa}_1^{i}(W_1, ..., W_i))=s^{-\frac{\epsilon_0}{8}(i-1)}cap(\Xi^A) \otimes_{k=1}^{i}P_{\bar{e}}[Y_{\cdot \wedge T_{\overline{\Xi^B}}} \in W_{k}].
\end{aligned}
\end{equation}
So that Lemma 7.2 in \cite{1} could also be generalized straightforwardly by Lemma \ref{L5.1}. Thus we have that there exist some $c, c(\epsilon_0)>0$ such that if $1\geq \delta \geq s^{-c(\epsilon_0)}$, $n\leq s^{c(\epsilon_0)}$, $u\geq s^{-c(\epsilon_0)}$ and $ s\geq c(\epsilon_0)$, then we can construct a probability space $(\Omega,\mathcal{A},Q)$ with independent Poisson point processes $\mu_2, \mu_3, ...$ and another Poisson point process $\theta$, where $\mu_i$ has intensity $u\kappa_1^i$ and $\theta$ has intensity $u\delta \kappa_2$, $u\geq 0$, such that
\begin{equation}
\label{6.60}
Q[\cup_{i\geq 2}\mathcal{T}(\mu_i) \subset \mathcal{T}(\theta)] \geq 1- ce^{-cu\delta cap(\Xi^A)}.
\end{equation}
Thus Proposition 4.2 in \cite{1} has been generalized into Lemma \ref{L3.2}. \qed \\
\section{Gumbel Fluctuations}
\indent In this section we will use Proposition \ref{T2} to prove Theorem \ref{T1}. The idea comes from Section 3 in \cite{1}.\\
\indent For convenience we introduce the partial uncovered subset $F_\rho$ by
\begin{equation}
\label{3.27}
F_{\rho}=F\setminus Y(0,t(\rho)), \text{ where } t(\rho)=N^d(1-\rho)g(0)\log|F| \text{ and } 0<\rho<1.
\end{equation}
\indent By Proposition \ref{T2}, we could generalize Lemma 3.1 and Lemma 3.2 in \cite{1} straightforwardly into the following two lemmas.
\begin{lemma}[a generalization of Lemma 3.1 in \cite{1}]
\label{L6.1}
\indent Given $d \geq 3$ and moving mode matrix $Q$, there exist constants $c_7>0$, $c_8>0$ and some $c>0$ such that if $F \subset \mathbb{T}_N$ satisfies $2 \leq |F| \leq N^{c_7}$
and $d_\infty(x,y)>|F|^{c_8}$ for all $x,y \in F, x\neq y$, then
\begin{equation}
\label{3.28}
\sup_{z\in \mathbb{R},x\in \mathbb{T}_N} |P_x[F\subset Y(0, u_F(z)N^d)]-e^{-e^{-z}}| \leq c|F|^{-c}.
\end{equation}
\end{lemma}
\begin{lemma}[a generalization of Lemma 3.2 in \cite{1}]
\label{L6.2}
\indent Given $d \geq 3$ and moving mode matrix $Q$, there exist constants $c_9>0$ and some $c, c(\rho)>0$ such that for $0<\rho \leq c_9$ and $F \subset \mathbb{T}_N$ with $|F| \geq c(\rho)$
\begin{equation}
\label{3.29}
P[F_{\rho} \notin \mathcal{G}_1] \leq c|F|^{-c(\rho)}.
\end{equation}
where $\mathcal{G}_1=\{ F' \subset F : ||F'|-|F|^{\rho}| \leq |F|^{\frac{2}{3}\rho}, \mathop{\inf}_{x,y\in F', x\neq y} d_\infty(x,y) \geq |F|^{\frac{1}{2d}} \}$.
\end{lemma}
\indent Now we give the proof of Theorem \ref{T1}.\\
\indent \emph{Proof of Theorem \ref{T1}} By (\ref{3.29}) we have for $0<\rho \leq c_9$ and $|F| \geq c(\rho)$,
\begin{equation}
\label{F1}
|P[T_C^F \leq u_F(z)N^d]-P[F\subset Y(0,u_F(z)N^d), F_{\rho}\in \mathcal{G}_1]|\leq c|F|^{-c(\rho)}.
\end{equation}
Also by the Markov property, if $|F| \geq c(\rho)$ (so that $\emptyset \notin \mathcal{G}_1$),
\begin{equation}
\label{F2}
\begin{aligned}
& P[F\subset Y(0,u_F(z)N^d), F_{\rho}\in \mathcal{G}_1]\\
= &\sum_{x\in \mathbb{T}_N, F'\in \mathcal{G}_1}P[F_\rho=F',Y_{t(\rho)}=x]P_x[F'\subset Y(0,u_F(z)N^d-t(\rho))].
\end{aligned}
\end{equation}
Set $h=\log \frac{|F'|}{|F|^{\rho}}$ so that $P_x[F' \subset Y(0,u_F(z)N^d-t(\rho))]=P_x[F'\subset Y(0,u_{F'}(z-h)N^d)]$. Also fix $\rho=c\leq c_9$ small enough so that $2d\rho \leq c_7$ and $\frac{1}{4d\rho} \geq c_8$. Then there exists $c>0$ such that (\ref{F1}) holds when $|F| \geq c$. Furthermore, we have that Lemma \ref{L6.1} applies to all $F' \in \mathcal{G}_1$ when $|F|\geq c$, since by the definition of $\mathcal{G}_1$ every $F' \in \mathcal{G}_1$ satisfies $|F'| \leq 2|F|^{\rho} \leq |F|^{2\rho}$, if $|F|\geq c=c(\rho)$. So that $2d\rho \leq c_7$ implies $|F'| \leq |F|^{2\rho} \leq N^{2d\rho} \leq N^{c_7}$ and $\frac{1}{4d\rho} \geq c_8$ implies $\inf_{x,y\in F', x\neq y}d_{\infty}(x,y) \geq |F|^{\frac{1}{2d}} \geq |F'|^{\frac{1}{4d\rho}} \geq |F'|^{c_8}$. Thus applying Lemma \ref{L6.1} with $F'$ in place of $F$, we get that for all $|F|\geq c, x\in \mathbb{T}_N$ and $F' \in \mathcal{G}_1$, there exist some $c>0$ such that
\begin{equation}
\label{F3}
|P_x[F'\subset Y(0,u_F(z)N^d-t(\rho))]-e^{-e^{-(z-h)}}|\leq c|F|^{-c}.
\end{equation}
But it is elementary that there exists a constant $c>0$ such that
\begin{equation}
|e^{-e^{-(z-h)}}-e^{-e^{-z}}| \leq c|h| \text{ for all } z,h \in \mathbb{R},
\end{equation}
and for $h$ we have there exist some $c>0$ such that $|h| \leq \max(\log(1+|F|^{-\frac{1}{3}\rho}),-\log(1-|F|^{-\frac{1}{3}\rho})) \leq c|F|^{-\frac{1}{3}\rho}$ provided $|F|>c$, so in fact there exist some $c>0$ such that $|P_x[F'\subset Y(0,u_F(z)N^d-t(\rho))]-e^{-e^{-z}}|\leq c|F|^{-c}$ for all $F' \in \mathcal{G}_1$. Thus (\ref{F2}) implies that there exist some $c>0$ such that
\begin{equation}
\label{F4}
|P[F\subset Y(0,u_F(z)N^d), F_\rho \in \mathcal{G}_1] -e^{-e^{-z}}P[F_\rho \in \mathcal{G}_1]| \leq c|F|^{-c}.
\end{equation}
Combining this with (\ref{F1}) and one more application of (\ref{3.29}), Theorem \ref{T1} holds when $|F| \geq c$ where $c>0$ is a constant. Furthermore, we could adjust constants to show that Theorem \ref{T1} holds for all $F$, so the proof of Theorem \ref{T1} is complete.\qed \\

\emph{Acknowledgement. The authors would like to thank Prof Xinyi Li for valuable suggestions.}

\end{document}